\documentclass[review]{elsarticle}
\usepackage{geometry}
\usepackage{lineno,hyperref}
\modulolinenumbers[5]

\usepackage{amssymb}
\usepackage[dvipsnames]{xcolor}
\usepackage{pgfplots}
\usepackage{tikz}
\usetikzlibrary{decorations.pathreplacing}
\usetikzlibrary{intersections}
\usepackage[fleqn]{amsmath}
\newtheorem{theorem}{Theorem}[section]

\newtheorem{lemma}[theorem]{Lemma}

\newtheorem{prop}[theorem]{Proposition}

\begin{document}

\begin{frontmatter}

%
%
%
%
%

\title{Bounds on some monotonic topological indices of bipartite\\
 graphs with a given number of cut edges\tnoteref{t1}}
\tnotetext[t1]{Project supported by
the National Natural Science Foundation of China (61572190) and
Hunan Provincial Innovation Foundation for Postgraduate.}

\author{Hanlin Chen}
\ead{hanlinchen@yeah.net}
\author{Renfang Wu}
\ead{wfang2005@126.com}
\author{Hanyuan Deng\corref{cor1}}
\ead{hydeng@hunnu.edu.cn}
\cortext[cor1]{Corresponding author}
\address{College of Mathematics and Computer Science, Hunan Normal University, Changsha, Hunan 410081, P. R. China}
\address{Key Laboratory of High Performance Computing and Stochastic Information Processing (HPCSIP) (Ministry of Education of China)}
\begin{abstract}
Let $I(G)$ be a topological index of a graph. If $I(G+e)<I(G)$ (or $I(G+e)>I(G)$, respectively) for each edge $e\not\in G$, then $I(G)$ is monotonically decreasing (or increasing, respectively) with the addition of edges. In this article, we present lower or upper bounds for some monotonic topological indices,
including the Wiener index, the hyper-Wiener index, the Harary index, the connective eccentricity index, the eccentricity distance sum of bipartite graphs in terms of the number of cut edges, and characterize the corresponding extremal graphs, respectively.
\end{abstract}


\end{frontmatter}

\linenumbers

\section{Introduction}

\label{intro}
Let $G$ be a simple connected graph with vertex set $V(G)$ and edge set $E(G)$.
For a vertex $u\in V(G)$,
denote by $N_G(u)$ the neighborhood of $u$ in $G$ and by $d_{G}(u) = |N_G(u)|$ the degree of $u$ in $G$.
A vertex of $G$ is called pendent if it has degree 1, and the edge incident
with a pendent vertex is a pendent edge.
 For vertices $u, v\in V(G)$, the distance $d_{G}(u,v)$ is defined as the length
of a shortest path between $u$ and $v$ in $G$. The eccentricity $\varepsilon_{G}(u)$ of a vertex $u$ is the maximum distance from $u$ to any other
vertex, i.e., $\varepsilon_{G}(u)=\max \{d_{G}(u,v)| v\in V(G)\}$. A cut edge is an edge whose deletion increases the number of connected
components. As usual, let $S_{n}$ and $K_{n}$ be a star and a complete graph with $n$ vertices, respectively.

A graph is bipartite if its vertex set can be partitioned into two subsets $X$ and $Y$ so that every edge has one end in $X$ and one end in $Y$; such a partition $(X, Y)$ is called a bipartition of the graph, and $X$ and $Y$ its parts. We denote a bipartite graph $G$ with bipartition $(X, Y)$ by $G[X, Y]$. If $G[X, Y]$ is simple and every vertex in $X$ is joined to every vertex in $Y$, then $G$ is called a complete bipartite graph, which is denoted by $K_{s,t}$, where $s=|X|$ and $t=|Y|$.

Topological indices are numerical parameters of a (molecular) graph which characterize its topology and are usually graph invariants. They have been used for
examining quantitative structure-activity relationships (QSARs) extensively in which the biological activity or other properties of molecules are correlated with their chemical structure. Thousands of topological indices have been developed describing structural, physicochemical properties or biological activity of molecular graphs.

It is well known that many important topological indices have the monotonicity \cite{chen1,chen2,wu1}, i.e., they decrease (or increase, respectively) with addition of edges, including the Wiener index, the Kirchhoff index, the Hosoya index, the Estrada index, the Zagreb index etc.

In recent years, many literatures have been involved in the research of the extremal values of some topological indices in terms of graph structure parameters. Among all the connected graphs with $n$ and $k$ cut edges, Hua \cite{hua2009} determined the minimal value of the Wiener index, and Xu and Trinajsti\'{c} \cite{xu2011} characterized the minimal graph with respect to the hyper-Wiener index and the maximal one with respect to the Harary index. Hua et al. \cite{hua2012} characterized the graphs with the minimum eccentricity distance sum within all connected graphs on $n$ vertices with $k$ cut edges and all graphs on $n$ vertices with edge-connectivity $k$, respectively. In fact, all these indices in these literatures above are monotonic. In light of the monotonicity, the present authors \cite{wu1} determined the extremal values of some monotonic topological indices in terms of the number of cut vertices, or number of cut edges, or the vertex connectivity, or the edge connectivity of a graph. In \cite{chen2}, the authors determined the extremal values of some monotonic topological indices in graphs with a given vertex bipartiteness. In \cite{chen1}, the authors determined the extremal values of some monotonic topological indices in all bipartite graphs with a given matching number.

Motivated from the above results, it is natural to consider these extremal problems from the class of general connected graphs to the bipartite graphs, and it is interesting to give a unified approach to these extremal problems. In this paper, we focus on some monotonic topological indices of bipartite graphs with a given number of cut edges.

\section{Basic properties}

We consider the effect edge addition (or deletion) on topological indices.

Let $I(G)$ be a topological index of a graph. If $I(G+e)<I(G)$ (or $I(G+e)>I(G)$, resp.) for each edge $e\not\in G$, then $I(G)$ is monotonically decreasing (or increasing, respectively) with the addition of edges.

This property is shared by many important topological indices. For example, the Wiener index, the Kirchhoff index, the eccentricity distance sum and the Merrifield-Simmons index are monotonically decreasing with the addition of edges, and the Zagreb index, the Hosoya index, the Estrada index and the ABC index are monotonically increasing with the addition of edges.

It obvious that the complete graph $K_{n}$ possess the maximum value for the monotonically increasing topological indices and the minimum value for the monotonically decreasing topological indices among all connected graphs with $n$ vertices.

\begin{prop}\label{p1}
Let $G$ be a bipartite graph with the minimal $I$-value for the monotonically decreasing topological index $I$ (or, the maximal $I$-value for the monotonically increasing topological index $I$) among all bipartite graphs with $n$ vertices and $k$ cut edges $e_{1},  e_{2}, \cdots, e_{k}$. Then each component of $G-\{e_{1}, e_{2},\cdots, e_{k}\}$ is a complete bipartite graph or a single vertex.
\end{prop}
\textbf{Proof}. Without loss of generality, we assume that $I$ is a monotonically decreasing topological index. Let $G$ be a connected bipartite graph with the minimal $I$-value among all the bipartite graphs with $n$ vertices and $k$ cut edges $e_{1},  e_{2}, \cdots, e_{k}$. If there is a non-trivial component $H$ (with at least four vertices) of $G-\{e_{1},  e_{2}, \cdots, e_{k}\}$ such that $H$ is not a complete bipartite graph, then, by adding an edge $e$ between a pair of vertices which come from different partitions in $H$, we can obtain a new bipartite graph $G+e$ with $n$ vertices and $k$ cut edges $e_{1},  e_{2}, \cdots, e_{k}$, and $I(G+e)<I(G)$, a contradiction.   \hfill $\Box$

Proposition \ref{p1} describes a common structural characteristic of the extremal graphs for monotonic topological indices over all bipartite graphs with $n$ vertices and $k$ cut edges.

\section{Bounds on some topological indices of bipartite graphs with $k$ cut edges}

In this section, we will give some bounds on the Wiener index, the hyper-Wiener index, the Harary index, the connective eccentricity index and the eccentricity distance sum of bipartite graphs with a given number of cut edges.

The Wiener index \cite{wiener1947} is one of the most used topological indices with high correlation with many physical and chemical of molecular compounds.
The Wiener index of a graph $G$, denoted by $W(G)$, is defined as the sum of all distances between any two vertices in $G$, that is
\begin{equation*}
  W(G)=\sum_{\{u,v\}\subseteq V(G)}d_{G}(u,v).
\end{equation*}

The hyper-Wiener index of a graph $G$, denoted by $WW(G)$, was first introduced by Randi\'{c} \cite{randic1993} and defined as
\begin{equation*}
  WW(G)=\frac{1}{2}\sum_{\{u,v\}\subseteq V(G)}[d_{G}(u,v)+d_{G}(u,v)^2].
\end{equation*}

The Harary index of a graph $G$, denoted by $H(G)$, has been introduced independently by Plav\u{s}i\'{c} et al. \cite{plavsic1993} and Ivanciuc et al. \cite{inanciuc1993} in 1993 for the characterization of molecular graphs. In contrast with the Wiener index, the Harary index is defined as
\begin{equation*}
  H(G)=\sum_{\{u,v\}\subseteq V(G)}\frac{1}{d_{G}(u,v)}.
\end{equation*}

The connective eccentricity index (CEI), as a novel, adjacency-cum-path length based, topological descriptor, was introduced by Gupta et al. \cite{gupta2000} in 2000, and is defined as
\begin{equation*}
  \xi^{ce}(G)=\sum_{uv\in E(G)}\bigg(\frac{1}{\varepsilon_{G}(u)}+\frac{1}{\varepsilon_{G}(v)}\bigg)
  =\sum_{u\in V(G)}\frac{d_{G}(u)}{\varepsilon_{G}(u)}.
\end{equation*}

The eccentricity distance sum (EDS) \cite{gupta2002} of a graph $G$ is defined as
\begin{equation*}
  \xi^{d}(G)=\sum_{\{u,v\}\subseteq V(G)}\big[\varepsilon_{G}(u)+\varepsilon_{G}(v)\big] d_{G}(u,v)
  =\sum_{u\in V(G)}\varepsilon_{G}(u) D_{G}(u),
\end{equation*}
where $D_{G}(u)=\sum_{v\in V(G)}d_{G}(u,v)$.

The following lemma is the direct consequence of the definitions above.

\begin{lemma}\label{l2}
Let $G$ be a connected graph of order $n$ and not isomorphic to complete graph $K_{n}$, and $e\notin E(G)$. Then \textup{(i)} $W(G+e)<W(G)$;
\textup{(ii)} $WW(G+e)<WW(G)$;
\textup{(iii)} $H(G+e)>H(G)$;  \textup{(iv)} $\xi^{ce}(G+e)>\xi^{ce}(G)$;  \textup{(v)} $\xi^{d}(G+e)<\xi^{d}(G)$.
\end{lemma}

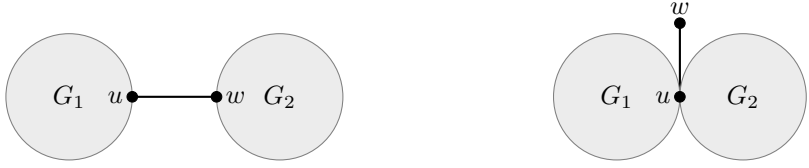
\begin{figure}[ht!]
\begin{center}
\begin{tikzpicture}[scale=1.4]
\draw [gray,fill=gray!15] (0,0) circle (.6cm);
\node  at (0,0) {$G_{1}$};
\draw [gray,fill=gray!15] (2,0) circle (.6cm);
\node  at (2,0) {$G_{2}$};
\draw [fill=black] (.6,0) circle (.05cm);
\draw [fill=black] (1.4,0) circle (.05cm);
\draw [thick]  (.6,0)--(1.4,0)(5.8,0)--(5.8,.7);
\draw [gray,fill=gray!15] (5.2,0) circle (.6cm);
\node  at (5.2,0) {$G_{1}$};
\draw [gray,fill=gray!15] (6.4,0) circle (.6cm);
\node  at (6.4,0) {$G_{2}$};
\draw [fill=black] (5.8,0) circle (.05cm);
\draw [fill=black] (5.8,.7) circle (.05cm);
\node  [left] at (.6,0) {$u$};
\node  [right] at (1.4,0) {$w$};
\node  [left] at (5.8,0) {$u$};
\node  [above] at (5.8,.7) {$w$};
\end{tikzpicture}
\caption{Graphs $G$ (\emph{left}) and $G'$ (\emph{right}) in Lemma \ref{l3}.}
\label{fig:1}
\end{center}
\end{figure}

\begin{lemma}\label{l3}
If $G$ is the graph obtained from $G_{1}\cup G_{2}$ by adding an edge $uw$, where $G_{1}$ and $G_{2}$ are two disjoint
connected graphs of order at least 2 with $u\in V(G_{1})$ and $w\in V(G_{2})$, and $G'$ is the graph obtained from $G_{1}\cup G_{2}$ by identifying $u$ and $w$ to a new vertex (say, $u$) and adding a pendent edge (say, $uw$ without confusion), see Figure 1. Then
\textup{(i)} \cite{hua2009} $W(G')<W(G)$; \textup{(ii)} \cite{xu2011} $WW(G')<WW(G)$;
 \textup{(iii)} \cite{xu2011} $H(G')>H(G)$;  \textup{(iv)} \cite{yu2013} $\xi^{ce}(G')>\xi^{ce}(G)$;  \textup{(v)} \cite{hua2012} $\xi^{d}(G')<\xi^{d}(G)$.
\end{lemma}

\begin{figure}[ht!]
\begin{center}
\begin{tikzpicture}[scale=1.3]
\draw [thick]  (0,0)--(0,1.8)--(1,0)--(1,1.8)--(2,0)--(2,1.8)--(3,0)--(3,1.8)--(2,0)--(0,1.8)--(3,0)--(1,1.8);
\draw [thick]  (1,1.8)--(0,0)--(2,1.8)(0,0)--(3,1.8)(2,1.8)--(1,0)--(3,1.8);
\draw [dashed] (1.8,2.3)--(2,1.8)--(2.2,2.3)(2,1.8)--(2,2.25);
\draw [dashed] (.8,2.3)--(1,1.8)--(1.2,2.3)(1,1.8)--(1,2.25);
\draw [dashed] (.8,-.5)--(1,0)--(1.2,-.5)(1,0)--(1,-.45);
\draw [dashed] (1.8,-.5)--(2,0)--(2.2,-.5)(2,0)--(2,-.45);
\draw [dashed] (2.8,-.5)--(3,0)--(3.2,-.5)(3,0)--(3,-.45);
\draw [dashed] (-.2,-.5)--(0,0)--(.2,-.5)(0,0)--(0,-.45);
\draw [fill=black] (0,1.8) circle (.05cm);
\node  at (.5,1.8) {$\cdots$};
\draw [fill=black] (1,1.8) circle (.05cm);
\node  at (1.5,1.8) {$\cdots$};
\draw [fill=black] (2,1.8) circle (.05cm);
\node  at (2.5,1.8) {$\cdots$};
\draw [fill=black] (3,1.8) circle (.05cm);
\draw [fill=black] (0,0) circle (.05cm);
\node  at (.5,0) {$\cdots$};
\draw [fill=black] (1,0) circle (.05cm);
\node  at (1.5,0) {$\cdots$};
\draw [fill=black] (2,0) circle (.05cm);
\node  at (2.5,0) {$\cdots$};
\draw [fill=black] (3,0) circle (.05cm);
\draw [fill=black] (5,1.8) circle (.05cm);
\node  at (5.5,1.8) {$\cdots$};
\draw [fill=black] (6,1.8) circle (.05cm);
\node  at (6.5,1.8) {$\cdots$};
\draw [fill=black] (7,1.8) circle (.05cm);
\node  at (7.5,1.8) {$\cdots$};
\draw [fill=black] (8,1.8) circle (.05cm);
\draw [fill=black] (5,0) circle (.05cm);
\node  at (5.5,0) {$\cdots$};
\draw [fill=black] (6,0) circle (.05cm);
\node  at (6.5,0) {$\cdots$};
\draw [fill=black] (7,0) circle (.05cm);
\node  at (7.5,0) {$\cdots$};
\draw [fill=black] (8,0) circle (.05cm);
\draw [thick]  (5,0)--(5,1.8)--(6,0)--(6,1.8)--(7,0)--(7,1.8)--(8,0)--(8,1.8)--(7,0)--(5,1.8)--(8,0)--(6,1.8);
\draw [thick]  (6,1.8)--(5,0)--(7,1.8)(5,0)--(8,1.8)(7,1.8)--(6,0)--(8,1.8);
\draw [thick]  (0.3,2.5)--(0,1.8)--(-.3,2.5)(0,1.8)--(-.6,2.5);
\draw [dashed] (6.8,2.3)--(7,1.8)--(7.2,2.3)(7,1.8)--(7,2.25);
\draw [dashed] (5.8,2.3)--(6,1.8)--(6.2,2.3)(6,1.8)--(6,2.25);
\draw [dashed] (5.8,-.5)--(6,0)--(6.2,-.5)(6,0)--(6,-.45);
\draw [dashed] (6.8,-.5)--(7,0)--(7.2,-.5)(7,0)--(7,-.45);
\draw [dashed] (7.8,-.5)--(8,0)--(8.2,-.5)(8,0)--(8,-.45);
\draw [dashed] (4.8,-.5)--(5,0)--(5.2,-.5)(5,0)--(5,-.45);
\draw [thick]  (3.6,2.5)--(3,1.8)--(3,2.5)(3,1.8)--(2.7,2.5);
\draw [thick]  (7.5,2.5)--(8,1.8)--(8.1,2.5)(8.4,2.5)--(8,1.8)--(9,2.5);
\draw [fill=black] (0.3,2.5) circle (.04cm);
\draw [fill=black] (-.3,2.5) circle (.04cm);
\node  at (0.02,2.5) {$\cdots$};
\draw [fill=black] (-.6,2.5) circle (.04cm);
\node [left]  at (0,1.8) {$u$};
\node [above]  at (-.65,2.5) {$u_{1}$};
\node [above]  at (-.3,2.5) {$u_{2}$};
\node [above]  at (.3,2.5) {$u_{a}$};
\node [right]  at (3,1.8) {$v$};
\draw [fill=black] (3.6,2.5) circle (.04cm);
\draw [fill=black] (3,2.5) circle (.04cm);
\node  at (3.32,2.5) {$\cdots$};
\draw [fill=black] (2.7,2.5) circle (.04cm);
\node [above]  at (2.65,2.5) {$v_{1}$};
\node [above]  at (3,2.5) {$v_{2}$};
\node [above]  at (3.6,2.5) {$v_{b}$};
\node [left]  at (5,1.8) {$u$};
\node [right]  at (8,1.8) {$v$};
\draw [fill=black] (8.4,2.5) circle (.04cm);
\node [above]  at (8.45,2.5) {$v_{1}$};
\node [above]  at (9.05,2.5) {$v_{b}$};
\draw [fill=black] (8.1,2.5) circle (.04cm);
\draw [fill=black] (7.5,2.5) circle (.04cm);
\node [above]  at (7.48,2.5) {$u_{1}$};
\node [above]  at (8.08,2.5) {$u_{a}$};
\node  at (7.82,2.5) {$\cdots$};
\node  at (8.72,2.5) {$\cdots$};
\draw [fill=black] (9,2.5) circle (.04cm);
\end{tikzpicture}
\caption{Graphs $G$ (\emph{left}) and $G'$ (\emph{right}) in Lemma \ref{l4}.}
\label{fig:2}
\end{center}
\end{figure}
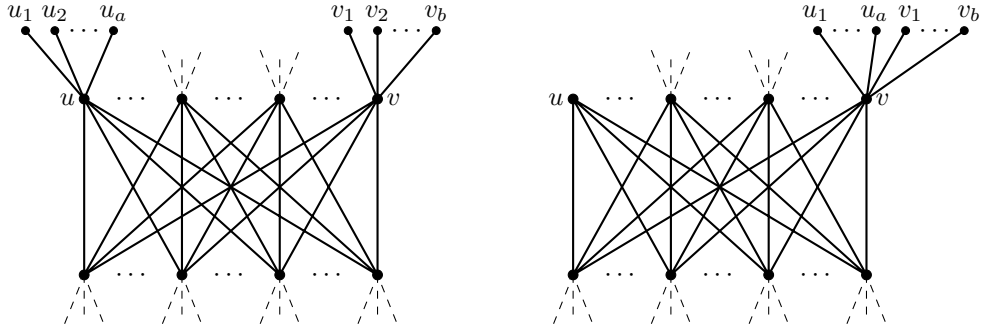

\begin{lemma}\label{l4}
 Let $G$ be a bipartite graph obtained from a complete bipartite graph $K_{s,t}$ by attaching some pendent vertices. If $u$ and $v$ are two vertices belonged to the same part of $K_{s,t}$, $a$ pendent vertices $u_{1}, u_{2},\cdots,u_{a}$ are attached to the vertex $u$ and $b$ pendent vertices $v_{1}, v_{2},\cdots,v_{b}$ are attached to the vertex $v$, $G'=G-\{uu_{1},uu_{2},\cdots,uu_{a}\}+\{vu_{1},vu_{2},\cdots,vu_{a}\}$, see Figure 2. Then
 \textup{(i)} $W(G')<W(G)$; \textup{(ii)} $WW(G')<WW(G)$;
\textup{(iii)} $H(G')>H(G)$;  \textup{(iv)} $\xi^{ce}(G')\geq\xi^{ce}(G)$;  \textup{(v)} $\xi^{d}(G')<\xi^{d}(G)$.
\end{lemma}
\textbf{Proof}. Firstly, for the Wiener index, we have
\begin{align*}
     W(G')-W(G)=&\sum_{i=1}^{a}\big[d_{G'}(u_{i},u)-d_{G}(u_{i},u)\big]+\sum_{i=1}^{a}\big[d_{G'}(u_{i},v)-d_{G}(u_{i},v)\big]
      +\sum_{i=1}^{a}\sum_{j=1}^{b}\big[d_{G'}(u_{i},v_{j})-d_{G}(u_{i},v_{j})\big]\\
     =&2a-2a-2ab\\
     =&-2ab<0.
\end{align*}

Similarly, one can obtain that $WW(G')-WW(G)=-7ab<0$ and $H(G')-H(G)=\frac{ab}{4}>0$.

Next, for the connective eccentricity index, we have the following two cases to consider.

{\bf Case I}. $\varepsilon_{G'}(u_{i})=\varepsilon_{G}(u_{i})-1$ for $i=1,2,\cdots,a$. Then no other vertices in the same part of $u$ and $v$ are attached by a pendent vertex.

Assume that $d_{G'}(u)=d$, then $d_{G}(u)=a+d$ and $d_{G'}(v)=d_{G}(v)+a=a+b+d$. Moreover, $\varepsilon_{G'}(u_{i})=\varepsilon_{G}(u_{i})-1=3$ for $i=1,2,\cdots,a$ and $\varepsilon_{G'}(v_{j})=\varepsilon_{G}(v_{j})-1=3$ for $j=1,2,\cdots,b$. We have
\begin{align*}
     \xi^{ce}(G')-\xi^{ce}(G) =&\frac{d_{G'}(u)}{\varepsilon_{G'}(u)}-\frac{d_{G}(u)}{\varepsilon_{G}(u)}
     +\frac{d_{G'}(v)}{\varepsilon_{G'}(v)}-\frac{d_{G}(v)}{\varepsilon_{G}(v)}
     +\sum_{i=1}^{a}\left(\frac{d_{G'}(u_{i})}{\varepsilon_{G'}(u_{i})}-\frac{d_{G}(u_{i})}{\varepsilon_{G}(u_{i})}\right)
      +\sum_{j=1}^{b}\left(\frac{d_{G'}(v_{j})}{\varepsilon_{G'}(v_{j})}-\frac{d_{G}(v_{j})}{\varepsilon_{G}(v_{j})}\right)\\
     =&\left(\frac{d}{3}-\frac{a+d}{3}\right)+\left(\frac{a+b+d}{3}-\frac{b+d}{3}\right)
     +\left(\frac{1}{3}-\frac{1}{4}\right)a+\left(\frac{1}{3}-\frac{1}{4}\right)b\\
     =&\frac{a+b}{12}>0.
\end{align*}

{\bf Case II}. $\varepsilon_{G'}(u_{i})=\varepsilon_{G}(u_{i})$ for $i=1,2,\cdots,a$. Then there are other vertices in the same part of $u$ and $v$ also attached by a pendent vertex, and $\varepsilon_{G'}(u_{i})=\varepsilon_{G}(u_{i})=4$ for $i=1,2,\cdots,a$.

In this case, the eccentricity of all the vertices are unchanged. Let $d_{G'}(u)=d$. Then $d_{G}(u)-a=d_{G'}(u)$, $d_{G}(v)+a=d_{G'}(v)$ and the degrees of other vertices remain unchanged. We have
\begin{align*}
     \xi^{ce}(G')-\xi^{ce}(G)=\frac{d_{G'}(u)}{\varepsilon_{G'}(u)}-\frac{d_{G}(u)}{\varepsilon_{G}(u)}
     +\frac{d_{G'}(v)}{\varepsilon_{G'}(v)}-\frac{d_{G}(v)}{\varepsilon_{G}(v)}
     =\frac{d}{3}-\frac{d+a}{3}+\frac{b+d}{3}-\frac{a+b+d}{3}=0.
\end{align*}
{\bf Remark}. For the graph $G$ above, if there are at least three vertices in the same part of $u$ and $v$ attached by a pendent vertex, then we can obtain a graph $G'$ such that there are exactly two vertices $u$ and $v$ attached by a pendent vertex by {\bf Case II} and $\xi^{ce}(G')=\xi^{ce}(G)$. Moreover, we can obtain another graph $G''$ such that $\xi^{ce}(G'')>\xi^{ce}(G)$ by {\bf Case I}. So, if there are at least three vertices in the same part of $G$, then we can always obtain a graph $G''$ such that $\xi^{ce}(G'')>\xi^{ce}(G)$.

Finally, for the eccentricity distance sum index, we also have the following two cases to consider.

{\bf Case I}.  $\varepsilon_{G'}(u_{i})=\varepsilon_{G}(u_{i})-1$ for $i=1,2,\cdots,a$.

Note that $\varepsilon_{G'}(u_{i})=\varepsilon_{G}(u_{i})-1=3$ for $i=1,2,\cdots,a$
and $\varepsilon_{G'}(v_{j})=\varepsilon_{G}(v_{j})-1=3$ for $j=1,2,\cdots,b$.  And it is easy to check that
$D_{G'}(u)=D_{G}(u)+2a$, $D_{G'}(v)=D_{G}(v)-2a$,
$D_{G'}(u_{i})=D_{G}(u_{i})-2b$ for $i=1,2,\cdots,a$ and $D_{G'}(v_{j})=D_{G}(v_{j})-2a$ for $j=1,2,\cdots,b$.
Then
\begin{align*}
     \xi^{d}(G')-\xi^{d}(G)=&\big[\varepsilon_{G'}(u)D_{G'}(u)-\varepsilon_{G}(u)D_{G}(u)\big]+
    \big[\varepsilon_{G'}(v)D_{G'}(v)-\varepsilon_{G}(v)D_{G}(v)\big]\\
     +&\sum_{i=1}^{a}\big[\varepsilon_{G'}(u_{i})D_{G'}(u_{i})-\varepsilon_{G}(u_{i})D_{G}(u_{i})\big]
     +\sum_{j=1}^{b}[\varepsilon_{G'}(v_{j})D_{G'}(v_{j})-\varepsilon_{G}(v_{j})D_{G}(v_{j})]\\
     =&6a-6a+\sum_{i=1}^{a}\big[3D_{G'}(u_{i})-4D_{G}(u_{i})\big]
     +\sum_{j=1}^{b}\big[3D_{G'}(v_{j})-4D_{G}(v_{j})\big]\\
     <&\sum_{i=1}^{a}\big[3D_{G'}(u_{i})-3D_{G}(u_{i})\big]
     +\sum_{j=1}^{b}\big[3D_{G'}(v_{j})-3D_{G}(v_{j})\big]\\
     =&-12ab<0.
\end{align*}

{\bf Case II}. $\varepsilon_{G'}(u_{i})=\varepsilon_{G}(u_{i})$ for $i=1,2,\cdots,a$.

In this case, the eccentricity of all the vertices are unchanged. Moreover, $D_{G'}(u)=D_{G}(u)+2a$, $D_{G'}(v)=D_{G}(v)-2a$,
$D_{G'}(u_{i})=D_{G}(u_{i})-2b$ for $i=1,2,\cdots,a$ and $D_{G'}(v_{j})=D_{G}(v_{j})-2a$ for $j=1,2,\cdots,b$. Then
\begin{align*}
     \xi^{d}(G')-\xi^{d}(G)=&[3D_{G'}(u)-3D_{G}(u)]+
    [3D_{G'}(v)-3D_{G}(v)]+\sum_{i=1}^{a}[4D_{G'}(u_{i})-4D_{G}(u_{i})]\\
     &+\sum_{j=1}^{b}[4D_{G'}(v_{j})-4D_{G}(v_{j})]\\
     =&6a-6a-8ab-8ab\\
     =&-16ab<0.
\end{align*}
This completes the proof.  \hfill $\Box$

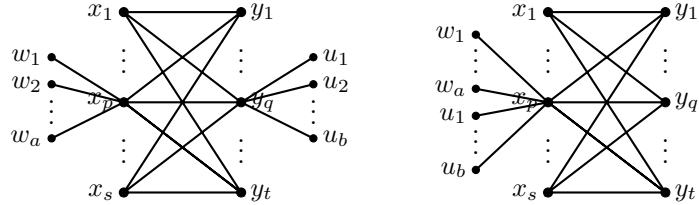
\begin{figure}[ht!]
\begin{center}
\begin{tikzpicture}[scale=1.2]
\draw [thick]  (0,0)--(1.3,0)--(0,1)--(1.3,1)--(0,2)--(1.3,2)--(0,1)--(1.3,0)--(0,2)(1.3,2)--(0,0)--(1.3,1);
\draw [thick]  (2.1,1.5)--(1.3,1)--(2.1,1.2)(1.3,1)--(2.1,.6);
\draw [thick]  (4.7,0)--(6,0)--(4.7,1)--(6,1)--(4.7,2)--(6,2)--(4.7,1)--(6,0)--(4.7,2)(6,2)--(4.7,0)--(6,1);
\draw [thick]  (3.9,1.75)--(4.7,1)--(3.9,1.15) (3.9,.85)--(4.7,1)--(3.9,.25);
\draw [fill=black] (0,0) circle (.05cm);
\node  [left] at (0,0) {$x_{s}$};
\node   at (0,.55) {$\vdots$};
\draw [fill=black] (0,1) circle (.05cm);
\node  [left] at (0,1) {$x_{p}$};
\node   at (0,1.55) {$\vdots$};
\draw [fill=black] (0,2) circle (.05cm);
\node  [left] at (0,2) {$x_{1}$};
\draw [fill=black] (1.3,0) circle (.05cm);
\node  [right] at (1.3,0) {$y_{t}$};
\node   at (1.3,.55) {$\vdots$};
\node  [right] at (1.3,1) {$y_{q}$};
\draw [fill=black] (1.3,1) circle (.05cm);
\node   at (1.3,1.55) {$\vdots$};
\draw [fill=black] (1.3,2) circle (.05cm);
\node  [right] at (1.3,2) {$y_{1}$};
\draw [fill=black] (-.8,1.5) circle (.04cm);
\node  [left] at (-.8,1.5) {$w_{1}$};
\draw [fill=black] (-.8,1.2) circle (.04cm);
\node  [left] at (-.8,1.2) {$w_{2}$};
\node   at (-.8,.98) {$\vdots$};
\draw [fill=black] (-.8,0.6) circle (.04cm);
\node  [left] at (-.8,0.6) {$w_{a}$};
\draw [thick]  (-.8,1.5)--(0,1)--(-.8,1.2)(-.8,0.6)--(0,1);
\draw [fill=black] (2.1,1.5) circle (.04cm);
\node  [right] at (2.1,1.5) {$u_{1}$};
\draw [fill=black] (2.1,1.2) circle (.04cm);
\node  [right] at (2.1,1.2) {$u_{2}$};
\node   at (2.1,.98) {$\vdots$};
\draw [fill=black] (2.1,0.6) circle (.04cm);
\node  [right] at (2.1,.6) {$u_{b}$};
\draw [fill=black] (4.7,0) circle (.05cm);
\node  [left] at (4.7,0) {$x_{s}$};
\node   at (4.7,.55) {$\vdots$};
\draw [fill=black] (4.7,1) circle (.05cm);
\node  [left] at (4.7,1) {$x_{p}$};
\node   at (4.7,1.55) {$\vdots$};
\draw [fill=black] (4.7,2) circle (.05cm);
\node  [left] at (4.7,2) {$x_{1}$};
\draw [fill=black] (6,0) circle (.05cm);
\node  [right] at (6,0) {$y_{t}$};
\node   at (6,.55) {$\vdots$};
\draw [fill=black] (6,1) circle (.05cm);
\node  [right] at (6,1) {$y_{q}$};
\node   at (6,1.55) {$\vdots$};
\draw [fill=black] (6,2) circle (.05cm);
\node  [right] at (6,2) {$y_{1}$};
\draw [fill=black] (3.9,1.75) circle (.04cm);
\node   at (3.9,1.55) {$\vdots$};
\node  [left] at (3.9,1.75) {$w_{1}$};
\draw [fill=black] (3.9,1.15) circle (.04cm);
\node  [left] at (3.9,1.15) {$w_{a}$};
\draw [fill=black] (3.9,.85) circle (.04cm);
\node   at (3.9,.65) {$\vdots$};
\node  [left] at (3.9,.85) {$u_{1}$};
\draw [fill=black] (3.9,.25) circle (.04cm);
\node  [left] at (3.9,.25) {$u_{b}$};
\end{tikzpicture}
\caption{ Graphs $G$ (\emph{left}) and $G'$ (\emph{right}) in Lemma \ref{l5}.}
\label{fig:3}
\end{center}
\end{figure}

\begin{lemma}\label{l5}
 Let $G$ be the bipartite graph arisen from $K_{s,t}$ by attaching $a$ pendent vertices to $x_{p}$ and $b$ pendent vertices to $y_{q}$, where the vertex $x_{p}$ belongs to the part with the cardinality $s$, the vertex $y_{q}$ belongs to the part with the cardinality $t$ of $K_{s,t}$  and $s\leq t$. $G'=G-\{y_{q}u_{1}, y_{q}u_{2},\cdots,y_{q}u_{b}\}+\{x_{p}u_{1}, x_{p}u_{2},\cdots,x_{p}u_{b}\}$, see Figure 3. Then
\textup{(i)} $W(G')<W(G)$;  \textup{(ii)} $WW(G')<WW(G)$;  \textup{(iii)} $H(G')>H(G)$;  \textup{(iv)} $\xi^{ce}(G')>\xi^{ce}(G)$;  \textup{(v)} $\xi^{d}(G')<\xi^{d}(G)$.
\end{lemma}
\textbf{Proof}. Firstly, by the definition of the Wiener index, we have
\begin{align*}
W(G')-W(G)=&\sum_{j=1}^{b}\sum_{i=1}^{a}[d_{G'}(u_{j},w_{i})-d_{G}(u_{j},w_{i})]
+\sum_{j=1}^{b}\sum_{k=1}^{s}[d_{G'}(u_{j},x_{k})-d_{G}(u_{j},x_{k})]\\
&+\sum_{j=1}^{b}\sum_{k=1}^{t}[d_{G'}(u_{j},y_{k})-d_{G}(u_{j},y_{k})]\\
=&-ab+[b(s-1)-1]+[1-b(t-1)]\\
<&b(s-t)\leq0.
\end{align*}
Analogously, it is easy to derive that $WW(G')-WW(G)<0$ and $H(G')-H(G)>0$.

Next, for the connective eccentricity index, we have
\begin{align*}
\xi^{ce}(G')-\xi^{ce}(G)=&\frac{d_{G'}(x_{p})}{\varepsilon_{G'}(x_{p})}-\frac{d_{G}(x_{p})}{\varepsilon_{G}(x_{p})}+
\frac{d_{G'}(y_{q})}{\varepsilon_{G'}(y_{q})}-\frac{d_{G}(y_{q})}{\varepsilon_{G}(y_{q})}
+\sum_{1\leq k\leq t,k\neq q}\bigg(\frac{d_{G'}(y_{k})}{\varepsilon_{G'}(y_{k})}-\frac{d_{G}(y_{k})}{\varepsilon_{G}(y_{k})}\bigg)\\
=&\frac{a+b+t}{2}-\frac{a+t}{2}
+\frac{s}{2}-\frac{b+s}{2}+\sum_{1\leq k\leq t,k\neq q}\Big(\frac{s}{2}-\frac{s}{3}\Big)\\
=&\frac{s(t-1)}{6}>0.
\end{align*}

Finally, it is easy to verify that
$D_{G'}(x_{p})=D_{G}(x_{p})-b$,  $D_{G'}(y_{q})=D_{G}(y_{q})+b$, $D_{G'}(w_{i})=D_{G}(w_{i})-b$ for $i=1,2,\cdots,a$ and $D_{G'}(u_{j})=D_{G}(u_{j})-a+s-t$ for $j=1,2,\cdots,b$. In addition,
$D_{G'}(x_{k})=D_{G}(x_{k})+b$ for $1\leq k\leq s$ and  $k\neq p$, $D_{G'}(y_{k})=D_{G}(y_{k})-b$ for $1\leq k\leq t$ and $k\neq q$.
Hence, this gives
\begin{align*}
\xi^{d}(G')-\xi^{d}(G)
=&\sum_{i=1}^{a}[3D_{G'}(w_{i})-3D_{G}(w_{i})]
+\sum_{j=1}^{b}[3D_{G'}(u_{j})-3D_{G}(u_{j})]+[2D_{G'}(x_{p})-2D_{G}(x_{p})]\\
+&\sum_{1\leq k\leq s,k\neq p}[3D_{G'}(x_{k})-3D_{G}(x_{k})]+[2D_{G'}(y_{q})-2D_{G}(y_{q})]
+\sum_{1\leq k\leq t,k\neq q}[2D_{G'}(y_{k})-3D_{G}(y_{k})]\\
<&-3ab+3b(s-t-a)-2b+3b(s-1)+2b-3b(t-1)\\
=&-6ab+6b(s-t)<0.
\end{align*}
This completes the proof. \hfill $\Box$

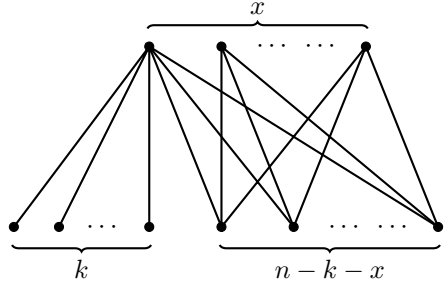
\begin{figure}[ht!]
\begin{center}
\begin{tikzpicture}[scale=1.2]
\draw [thick]  (0,-2)--(0,0)--(-1,-2)(-1.5,-2)--(0,0)--(.8,-2)--(.8,0)--(1.6,-2)--(2.4,0)--(3.2,-2)--(.8,0);
\draw [thick]  (.8,-2)--(2.4,0) (1.6,-2)--(0,0)--(3.2,-2);
\draw [fill=black] (0,0) circle (.05cm);
\draw [fill=black] (.8,0) circle (.05cm);
\draw [fill=black] (2.4,0) circle (.05cm);
\node  at (1.65,0) {$\cdots~\cdots$};
\draw [fill=black] (.8,-2) circle (.05cm);
\draw [fill=black] (1.6,-2) circle (.05cm);
\draw [fill=black] (3.2,-2) circle (.05cm);
\node  at (2.45,-2) {$\cdots~\cdots$};
\draw [fill=black] (0,-2) circle (.05cm);
\draw [fill=black] (-1,-2) circle (.05cm);
\draw [fill=black] (-1.5,-2) circle (.05cm);
\node  at (-.5,-2) {$\cdots$};
\draw[thick,decorate,decoration=brace] (-0.02,0.2) --node[above=0.4ex] {$x$} (2.42,0.2);
\draw[thick,decorate,decoration=brace] (0.02,-2.2) --node[below=0.4ex] {$k$} (-1.52,-2.2);
\draw[thick,decorate,decoration=brace] (3.22,-2.2) --node[below=0.4ex] {$n-k-x$} (.78,-2.2);
\end{tikzpicture}
\caption{The bipartite graph $B_{k}(x,n-k-x)$.}
\label{fig:4}
\end{center}
\end{figure}

Let $B_{k}(x,n-k-x)$ denote the bipartite graph arisen by attaching $k$ pendent vertices to a vertex with degree $n-k-x$ in $K_{x,n-k-x}$, where $x$, $k$ and $n$ are positive integers and $2\leq x\leq n-k-x$, see Figure 4.

\begin{theorem}\label{t6}
Let $G$ be a connected bipartite graph with $n\geq5$ vertices and $k$ cut edges.
\begin{enumerate}[(\romannumeral1)]
  \item If $k=n-1$, then $W(G)=\frac{n(n-1)}{2}$ and the graph $G$ is $S_{n}$;
  \item If $\frac{n-4}{2}\leq k\leq n-4$, then
 $W(G)\geq n^2-3n+2k+4$ with equality if and only if $G\cong B_{k}(2, n-k-2)$;
  \item If $1\leq k<\frac{n-4}{2}$ and $n$ is odd, then
$W(G)\geq\frac{3n^2+1}{4}+kn-k^2-2k-n$ with
equality if and only if $G\cong B_{k}(\frac{n-2k-1}{2}, \frac{n+1}{2})$ or
$G\cong B_{k}(\frac{n-2k+1}{2}, \frac{n-1}{2})$;
  \item If $1\leq k<\frac{n-4}{2}$ and $n$ is even, then $W(G)\geq
\frac{3n^2}{4}+kn-k^2-2k-n$
with equality if and only if $G\cong B_{k}(\frac{n-2k}{2}, \frac{n}{2})$.
\end{enumerate}
\end{theorem}
\textbf{Proof}. If $k=n-1$, the unique bipartite graph is $S_{n}$, and $W(S_{n})=\frac{n(n-1)}{2}$, as desired.

Note that $k\neq n-2$ and $k\neq n-3$ for any bipartite graph. In the following, we only need to consider $1\leq k\leq n-4$.
Let $G$ be a bipartite graph with the minimum Wiener index among all bipartite graphs with $n$ vertices and $k$ cut edges $e_{1},e_{2},\cdots,e_{k}$.
By Lemma 2 and Proposition 1, each component of $G-\{e_{1}, e_{2},\cdots,e_{k}\}$ is a complete bipartite graph or a single vertex.
By Lemma 3, all of its cut edges must be pendent edges. It is evident that the vertices, except all its pendent vertices,
form a complete bipartite graph by using Lemma 2 again. Moreover, we can confirm that all these pendent edges are
attached to the same vertex (i.e., $G\cong B_{k}(x,n-k-x)$, where $2\leq x\leq n-k-x$) from Lemma 4 and Lemma 5. And
\begin{align*}
  W(B_{k}(x,n-k-x))=f(x)=x^2+(2k-n)x+n^2-n-2k.
\end{align*}

If $\frac{n-4}{2}\leq k\leq n-4$, then by the fact $2\leq x\leq n-k-x$, we have
\begin{equation*}
 W(G)\geq f(2)=n^2-3n+2k+4
\end{equation*}
with equality if and only if $G\cong B_{k}(2, n-k-2)$.

If $1\leq k<\frac{n-4}{2}$, then
\begin{equation*}
  \min f(x)=
\begin{cases}
 f(\frac{n-2k-1}{2})=f(\frac{n-2k+1}{2}), & \mbox{if $n$ is odd};\\
f(\frac{n-2k}{2}), & \mbox{if $n$ is even}.
\end{cases}
\end{equation*}
Therefore, we have (1) $W(G)\geq f(\frac{n-2k-1}{2})=f(\frac{n-2k+1}{2})=\frac{3n^2+1}{4}+kn-k^2-2k-n$ for odd $n$,
with equality if and only if $G\cong B_{k}(\frac{n-2k-1}{2}, \frac{n+1}{2})$ or $G\cong B_{k}(\frac{n-2k+1}{2}, \frac{n-1}{2})$; and (2) $W(G)\geq f(\frac{n-2k}{2})=\frac{3n^2}{4}+kn-k^2-2k-n$ for even $n$, with equality if and only if $G\cong B_{k}(\frac{n-2k}{2}, \frac{n}{2})$.

This completes the proof. \hfill $\Box$

\begin{theorem}\label{t7}
Let $G$ be a connected bipartite graph with $n\geq5$ vertices and $k$ cut edges.
\begin{enumerate}[(\romannumeral1)]
  \item If $k=n-1$, then $WW(G)=\frac{(n-1)(3n-4)}{2}$ and the graph $G$ is $S_{n}$;
  \item If $\frac{2n-8}{5}\leq k\leq n-4$, then
 $WW(G)\leq \frac{3n^2-11n}{2}+5k+8$ with equality if and only if $G\cong B_{k}(2, n-k-2)$;
  \item If $1\leq k<\frac{2n-8}{5}$ and $n-5k\equiv 0~(mod ~4)$, then
$WW(G)\leq
n^2+\frac{5}{2}kn-\frac{3}{2}n-5k-\frac{25k^2}{8}$
with equality if and only if $G\cong B_{k}(\frac{2n-5k}{4}, \frac{2n+k}{4})$;
  \item If $1\leq k<\frac{2n-8}{5}$ and $n-5k\equiv 1~(mod~ 4)$, then $WW(G)\leq
n^2+\frac{5}{2}kn-\frac{3}{2}n-5k-\frac{25k^2-1}{8}$
with equality if and only if
$G\cong B_{k}(\frac{2n-5k-1}{4}, \frac{2n+k+1}{4})$;
  \item If $1\leq k<\frac{2n-8}{5}$ and $n-5k\equiv 2~(mod~4)$, then $WW(G)\leq
n^2+\frac{5}{2}kn-\frac{3}{2}n-5k-\frac{25k^2-4}{8}$
with equality if and only if $G\cong B_{k}(\frac{2n-5k-2}{4}, \frac{2n+k+2}{4})$
or $G\cong B_{k}(\frac{2n-5k+2}{4}, \frac{2n+k-2}{4})$;
\item If $1\leq k<\frac{2n-8}{5}$ and $n-5k\equiv 3~(mod~ 4)$, then $WW(G)\leq
n^2+\frac{5}{2}kn-\frac{3}{2}n-5k-\frac{25k^2-1}{8}$
with equality if and only if
$G\cong B_{k}(\frac{2n-5k+1}{4}, \frac{2n+k-1}{4})$.
\end{enumerate}
\end{theorem}
\textbf{Proof.}
If $k=n-1$, then the unique bipartite graph is $S_{n}$, and $WW(S_{n})=\frac{(n-1)(3n-4)}{2}$, as desired.

Let $G$ be a bipartite graph with the minimum hyper-Wiener index among all bipartite graphs with $n$ vertices and $k$ cut edges. By the same way as in the proof of Theorem \ref{t6}, we can confirm that $G\cong B_{k}(x,n-k-x)$, and
\begin{align*}
WW(B_{k}(x,n-k-x))=g(x)=2x^2+(5k-2n)x+\frac{3}{2}n^2-\frac{3}{2}n-5k.
\end{align*}
Also, it is easy to see that
$\min g(x)=g(2)$ for $\frac{2n-8}{5}\leq k\leq n-4$, i.e., if $\frac{2n-8}{5}\leq k\leq n-4$, then
\begin{equation*}
 WW(G)\geq g(2)=\frac{3n^2-11n}{2}+5k+8
\end{equation*}
with equality if and only if $G\cong B_{k}(2, n-k-2)$.

If $1\leq k<\frac{2n-8}{5}$, then
\begin{equation*}
  \min g(x)=
  \begin{cases}
g(\frac{2n-5k}{4}), & \mbox{if $n-5k\equiv 0$ (mod 4)};\\
g(\frac{2n-5k-1}{4}), & \mbox{if $n-5k\equiv 1$ (mod 4)};\\
g(\frac{2n-5k-2}{4})=g(\frac{2n-5k+2}{4}), & \mbox{if $n-5k\equiv 2$ (mod 4)};\\
g(\frac{2n-5k+1}{4}), & \mbox{if $n-5k\equiv 3$ (mod 4)}.
\end{cases}
\end{equation*}
Therefore, we have $WW(G)\geq g(\frac{2n-5k}{4})=n^2+\frac{5}{2}kn-\frac{3}{2}n-5k-\frac{25k^2}{8}$ for $n-5k\equiv 0~(mod ~4)$, with equality if and only if $G\cong B_{k}(\frac{2n-5k}{4}, \frac{2n+k}{4})$. Similarly, we can obtain $(iv)-(vi)$ for $n-5k\equiv 1~(mod ~4)$, $n-5k\equiv 2~(mod ~4)$ and
$n-5k\equiv 3~(mod ~4)$, respectively. \hfill $\Box$

\begin{theorem}\label{t8}
Let $G$ be a connected bipartite graph with $n\geq5$ vertices and $k$ cut edges.
\begin{enumerate}[(\romannumeral1)]
  \item If $k=n-1$, then $H(G)=\frac{n^2+n-2}{4}$ and the graph is $S_{n}$;
  \item If $\frac{3n-12}{4}\leq k\leq n-4$, then
 $H(G)\leq \frac{3n^2+9n-8k-24}{12}$ with equality if and only if $G\cong B_{k}(2, n-k-2)$;
  \item If $1\leq k<\frac{3n-12}{4}$ and $3n-4k\equiv 0 ~(mod~6)$, then
$H(G)\leq
\frac{3n^2-2n}{8}+\frac{2k^2-3kn+6k}{9}$
with equality if and only if $G\cong B_{k}(\frac{3n-4k}{6}, \frac{3n-2k}{6})$;
  \item If $1\leq k<\frac{3n-12}{4}$ and $3n-4k\equiv 1 ~(mod~6)$, then $H(G)\leq
\frac{3n^2-2n}{8}+\frac{2k^2-3kn+6k}{9}-\frac{1}{72}$ with
 equality if and only if $G\cong B_{k}(\frac{3n-4k-1}{6}, \frac{3n-2k+1}{6})$;
  \item If $1\leq k<\frac{3n-12}{4}$ and $3n-4k\equiv 2 ~(mod~6)$, then $H(G)\leq
\frac{3n^2-2n}{8}+\frac{2k^2-3kn+6k}{9}-\frac{1}{18}$ with
 equality if and only if $G\cong B_{k}(\frac{3n-4k-2}{6}, \frac{3n-2k+2}{6})$;
  \item If $1\leq k<\frac{3n-12}{4}$ and $3n-4k\equiv 3 ~(mod~6)$, then $H(G)\leq
\frac{3n^2-2n}{8}+\frac{2k^2-3kn+6k}{9}-\frac{1}{8}$ with
 equality if and only if $G\cong B_{k}(\frac{3n-4k-3}{6}, \frac{3n-2k+3}{6})$
or $G\cong B_{k}(\frac{3n-4k+3}{6}, \frac{3n-2k-3}{6})$;
\item  If $1\leq k<\frac{3n-12}{4}$ and $3n-4k\equiv 4 ~(mod~6)$, then $H(G)\leq
\frac{3n^2-2n}{8}+\frac{2k^2-3kn+6k}{9}-\frac{1}{18}$ with
 equality if and only if $G\cong B_{k}(\frac{3n-4k+2}{6}, \frac{3n-2k-2}{6})$;
  \item  If $1\leq k<\frac{3n-12}{4}$ and $3n-4k\equiv 5 ~(mod~6)$, then $H(G)\leq
\frac{3n^2-2n}{8}+\frac{2k^2-3kn+6k}{9}-\frac{1}{72}$ with
equality if and only if
$G\cong B_{k}(\frac{3n-4k+1}{6}, \frac{3n-2k-1}{6})$.
\end{enumerate}
\end{theorem}
\textbf{Proof}. If $k=n-1$, then the unique bipartite graph $G$ is $S_{n}$, and $H(S_{n})=\frac{n^2+n-2}{4}$, as desired.

Let $G$ be a bipartite graph with the maximum Harary index among all bipartite graphs with $n\geq5$ vertices and $k$ cut edges. By the same way as in the proof of Theorem \ref{t6}, we can confirm that $G\cong B_{k}(x,n-k-x)$, and
\begin{align*}
H(B_{k}(x,n-k-x))=h(x)=\frac{-6x^2+(6n-8k)x+3n^2-3n+8k}{12}.
\end{align*}

For $\frac{3n-12}{4}\leq k\leq n-4$, we can obtain that $\max h(x)=h(2)$, i.e.,
\begin{equation*}
  H(G)\leq h(2)=\frac{3n^2+9n-8k-24}{12}
\end{equation*}
with equality if and only if $G\cong B_{k}(2, n-k-2)$.

For $1\leq k<\frac{3n-12}{4}$, it is not difficulty to verify that
\begin{equation*}
  \max h(x)=
\begin{cases}
h(\frac{3n-4k}{6}), & \mbox{if $3n-4k\equiv 0$ (mod 6)};\\
h(\frac{3n-4k-1}{6}), & \mbox{if $3n-4k\equiv 1$ (mod 6)};\\
h(\frac{3n-4k-2}{6}), & \mbox{if $3n-4k\equiv 2$ (mod 6)};\\
h(\frac{3n-4k-3}{6})=h(\frac{3n-4k+3}{6}), & \mbox{if $3n-4k\equiv 3$ (mod 6)};\\
h(\frac{3n-4k+2}{6}), & \mbox{if $3n-4k\equiv 4$ (mod 6)};\\
h(\frac{3n-4k+1}{6}), & \mbox{if $3n-4k\equiv 5$ (mod 6)}.
\end{cases}
\end{equation*}
Hence, we have $H(G)\leq h(\frac{3n-4k}{6})=\frac{3n^2-2n}{8}+\frac{2k^2-3kn+6k}{9}$ for $3n-4k\equiv 0~ (mod~ 6)$,
with equality holds if and only if $G\cong B_{k}(\frac{3n-4k}{6}, \frac{3n-2k}{6})$. In a similar way, one can obtain $(iv)-(viii)$, respectively. \hfill $\Box$

\begin{theorem}\label{t9}
Let $G$ be a connected bipartite graph with $n\geq5$ vertices and $k$ cut edges.
\begin{enumerate}[(\romannumeral1)]
  \item  If $k=n-1$, then $\xi^{ce}(G)=\frac{n^2+n-2}{4}$ and the graph $G$ is $S_{n}$;
  \item If $1\leq k\leq n-4$ and $5n-5k\equiv 0~(mod~10)$, then $\xi^{ce}(G)\leq
\frac{5n^2+5k^2}{24}+\frac{n-5kn}{12}+\frac{3k}{4}$
with equality if and only if $G\cong B_{k}(\frac{5n-5k}{10}, \frac{5n-5k}{10})$;
  \item If $1\leq k\leq n-4$ and $5n-5k\equiv 5~(mod~10)$, then $\xi^{ce}(G)\leq
\frac{5n^2+5k^2}{24}+\frac{n-5kn}{12}+\frac{3k}{4}-\frac{1}{8}$
with equality  if and only if
$G\cong B_{k}(\frac{5n-5k-5}{10}, \frac{5n-5k+5}{10})$
\end{enumerate}
\end{theorem}
\textbf{Proof}. If $k=n-1$, then the unique bipartite graph $G$ is $S_{n}$, and $\xi^{ce}(G)=\frac{n^2+n-2}{4}$,
as desired.

Let $G$ be a bipartite graph with the maximum connective eccentricity index among all bipartite graphs with $n$ vertices and $k$ cut edges. By the same way as in the proof of Theorem \ref{t6}, we can confirm that $G\cong B_{k}(x,n-k-x)$, and
\begin{align*}
\xi^{ce}(B_{k}(x,n-k-x))=\psi(x)=\frac{-5x^2+(5n-5k-1)x+4k+n}{6}.
\end{align*}

For $1\leq k\leq n-4$, it is easy to verify that
\begin{equation*}
  \max \psi(x)=
\begin{cases}
\psi(\frac{5n-5k}{10}), & \mbox{if $5n-5k\equiv 0$ (mod 10)};\\
\psi(\frac{5n-5k-5}{10}), & \mbox{if $5n-5k\equiv 5$ (mod 10)}.
\end{cases}
\end{equation*}
Therefore, $\xi^{ce}(G)\leq \psi(\frac{5n-5k}{10})=
\frac{5n^2+5k^2}{24}+\frac{n-5kn}{12}+\frac{3k}{4}$ for $5n-5k\equiv 0~(mod~10)$, with equality if and only if $G\cong B_{k}(\frac{5n-5k}{10}, \frac{5n-5k}{10})$.
If $1\leq k<n-4$ and $5n-5k\equiv 5~(mod~10)$, then
 $\xi^{ce}(G)\leq \psi(\frac{5n-5k-5}{10})=\frac{5n^2+5k^2}{24}+\frac{n-5kn}{12}+\frac{3k}{4}-\frac{1}{8}$,
with equality if and only if $G\cong B_{k}(\frac{5n-5k-5}{10}, \frac{5n-5k+5}{10})$. \hfill $\Box$

\begin{theorem}\label{t10}
Let $G$ be a connected bipartite graph with $n\geq5$ vertices and $k$ cut edges.
\begin{enumerate}[(\romannumeral1)]
\item If $k=n-1$, then $\xi^{d}(G)=\frac{n^2+n-2}{4}$ and the graph $G$ is $S_{n}$;
\item If $\frac{3n-23}{11}\leq k\leq n-4$, then
 $\xi^{d}(G)\geq 4n^2+2kn-11n+8k+16$, with equality if and only if $G\cong B_{k}(2, n-k-2)$;
\item If $1\leq k<\frac{3n-23}{11}$ and $3n-11k+3\equiv 0~(mod~10)$, then $\xi^{d}(G)\geq\frac{71n^2-121k^2+31}{20}+\frac{53kn-59n-107k}{10}$
with equality if and only if $G\cong B_{k}(\frac{3n-11k+3}{10}, \frac{7n+k-3}{10})$;
\item If $1\leq k<\frac{3n-23}{11}$ and $3n-11k+3\equiv 1~(mod~10)$, then $\xi^{d}(G)\geq\frac{71n^2-121k^2+32}{20}+\frac{53kn-59n-107k}{10}$
with equality if and only if $G\cong B_{k}(\frac{3n-11k+2}{10}, \frac{7n+k-2}{10})$;
\item If $1\leq k<\frac{3n-23}{11}$ and $3n-11k+3\equiv 2~(mod~10)$, then $\xi^{d}(G)\geq\frac{71n^2-121k^2+35}{20}+\frac{53kn-59n-107k}{10}$
with equality if and only if $G\cong B_{k}(\frac{3n-11k+1}{10}, \frac{7n+k-1}{10})$;
\item If $1\leq k<\frac{3n-23}{11}$ and $3n-11k+3\equiv 3~(mod~10)$, then $\xi^{d}(G)\geq\frac{71n^2-121k^2+40}{20}+\frac{53kn-59n-107k}{10}$
with equality if and only if $G\cong B_{k}(\frac{3n-11k}{10}, \frac{7n+k}{10})$;
\item If $1\leq k<\frac{3n-23}{11}$ and $3n-11k+3\equiv 4~(mod~10)$, then $\xi^{d}(G)\geq\frac{71n^2-121k^2+47}{20}+\frac{53kn-59n-107k}{10}$
with equality if and only if $G\cong B_{k}(\frac{3n-11k-1}{10}, \frac{7n+k+1}{10})$;
\item If $1\leq k<\frac{3n-23}{11}$ and $3n-11k+3\equiv 5~(mod~10)$, then $\xi^{d}(G)\geq\frac{71n^2-121k^2+56}{20}+\frac{53kn-59n-107k}{10}$
with equality if and only if $G\cong B_{k}(\frac{3n-11k-2}{10}, \frac{7n+k+2}{10})$
or $G\cong B_{k}(\frac{3n-11k+8}{10}, \frac{7n+k-8}{10})$;
\item If $1\leq k<\frac{3n-23}{11}$ and $3n-11k+3\equiv 6~(mod~10)$, then $\xi^{d}(G)\geq\frac{71n^2-121k^2+47}{20}+\frac{53kn-59n-107k}{10}$
with equality if and only if $G\cong B_{k}(\frac{3n-11k+7}{10}, \frac{7n+k-7}{10})$;
\item If $1\leq k<\frac{3n-23}{11}$ and $3n-11k+3\equiv 7~(mod~10)$ , then $\xi^{d}(G)\geq\frac{71n^2-121k^2+40}{20}+\frac{53kn-59n-107k}{10}$
with equality if and only if $G\cong B_{k}(\frac{3n-11k+6}{10}, \frac{7n+k-6}{10})$;
\item If $1\leq k<\frac{3n-23}{11}$ and $3n-11k+3\equiv 8~(mod~10)$, then $\xi^{d}(G)\geq\frac{71n^2-121k^2+35}{20}+\frac{53kn-59n-107k}{10}$
with equality if and only if $G\cong B_{k}(\frac{3n-11k+5}{10}, \frac{7n+k-5}{10})$;
\item If $1\leq k<\frac{3n-23}{11}$ and $3n-11k+3\equiv 9~(mod~10)$, then $\xi^{d}(G)\geq\frac{71n^2-121k^2+32}{20}+\frac{53kn-59n-107k}{10}$
with equality if and only if $G\cong B_{k}(\frac{3n-11k+4}{10}, \frac{7n+k-4}{10})$.
\end{enumerate}
\end{theorem}
\textbf{Proof}. If $k=n-1$, then the unique bipartite graph is $S_{n}$, and $\xi^{d}(S_{n})=\frac{n^2+n-2}{4}$, as desired.

Let $G$ be a bipartite graph with the minimum eccentricity distance sum among all bipartite graphs with $n$ vertices and $k$ cut edges. By the same way as in the proof of Theorem \ref{t6}, we can confirm that $G\cong B_{k}(x,n-k-x)$, and
\begin{equation*}
\xi^{d}(B_{k}(x,n-k-x))=\varphi(x)=5x^2+(11k-3n-3)x+4n^2+2kn-5n-14k+2.
\end{equation*}

Clearly, for $\frac{3n-23}{11}\leq k\leq n-4$, $\varphi(x)\geq \varphi(2)$, that is
\begin{equation*}
  \xi^{d}(G)\geq \varphi(2)=4n^2+2kn-11n+8k+16
\end{equation*}
with equality if and only if $G\cong B_{k}(2, n-k-2)$.

If $1\leq k<\frac{3n-23}{11}$, then $\min \varphi(x)=\varphi(\frac{3n-11k+3}{10})$ for $3n-11k+3\equiv 0~ (mod ~10)$,
$\min \varphi(x)=\varphi(\frac{3n-11k+2}{10})$ for $3n-11k+3\equiv 1~ (mod ~10)$, etc. And (iii)-(xii) can be obtained, respectively, we omit them here.
\hfill $\Box$

\section{Conclusions}

In this paper, we first present a common structural characteristic of the extremal graphs for monotonic topological indices over all bipartite graphs with $n$ vertices and $k$ cut edges, and then determine the lower or upper bounds on the Wiener index, the hyper-Wiener index, the Harary index, the connective eccentricity index and the eccentricity distance sum of all bipartite graphs with a given number of cut edges, and characterize the corresponding extremal graphs.  Moreover, the methods we propose in this paper can be extended to determine the extremal values of other monotonic topological indices in bipartite graphs with a fixed number of cut edges. Along this line, some other interesting extremal problems on bipartite graphs with given parameters are valuable to be considered.





\end{document}